\documentclass{article}
\usepackage{amssymb}
\usepackage{latexsym}
\usepackage{amsfonts}
\usepackage{amsthm}
\usepackage{amsmath}

\usepackage{amssymb, amscd}

\newtheorem{problem}{Problem}[section]
\newtheorem{theo}[problem]{Theorem}

\newtheorem{defin}[problem]{Definition}

\newtheorem{cor}[problem]{Corollary}
\newtheorem{lema}[problem]{Lemma}

\begin{document}
\date{}
\title{{\LARGE Symmetric products of surfaces;\\
a unifying theme for topology and physics }}

\author{{\Large Pavle Blagojevi\' c} \\
        {Mathematics Institute SANU, Belgrade} \\[2mm]
{\Large Vladimir Gruji\' c} \\ {Faculty of Mathematics,
Belgrade}\\[2mm]
  {\Large  Rade \v Zivaljevi\' c}
 \\ {Mathematics Institute
 SANU, Belgrade}}
\maketitle

\begin{abstract}
This is a review paper about symmetric products of spaces
$SP^n(X):= X^n/S_n$. We focus our attention on the case of
$2$-manifolds $X$ and make a journey through selected topics of
algebraic topology, algebraic geometry, mathematical physics,
theoretical mechanics  etc. where these objects play an important
role, demonstrating along the way the fundamental unity of diverse
fields of physics and mathematics.
\end{abstract}

\section{Introduction}
\label{intro}

In recent years we have all witnessed a remarkable and extremely
stimulating exchange of deep and sophisticated ideas  between
Geometry and Physics and in particular between quantum physics and
topology. The student or a young scientist in one of these fields
is often urged to master elements of the other field as quickly as
possible, and to develop basic skills and intuition necessary for
understanding the contemporary research in both disciplines. The
topology and geometry of manifolds plays a central role in
mathematics and likewise in physics. The understanding of duality
phenomena for manifolds, mastering the calculus of characteristic
classes, as well as understanding the role of fundamental
invariants like the signature or Euler characteristic are just
examples of what is on the beginning of the growing list of
prerequisites for a student in these areas.

A graduate student of mathematics alone is often in position to
take many specialized courses covering different aspects of
manifold theory and related areas before an unified picture
emerges and she or he reaches the necessary level of maturity.

Obviously this state of affairs is somewhat unsatisfactory and
clearly it can only get worse in the future. One of objectives of
our review is to follow a single mathematical object, a remarkable
$(2n)$-dimensional manifold $SP^n(M)$, the $n$-th {\em symmetric
product} of a surface $M$, on a guided tour ``transversely''
through mathematics with occasional contacts with physics. We hope
that the reader will find this trip amusing and the information
interesting and complementary to the usual textbook presentations.

By definition, the $n$-th symmetric product of a space $X$ is
defined as $SP^n(X) := X^n/S_n$. In other words a point $D\in
SP^n(X)$ is an unordered collection of $n$-points in $X$, often
denoted by $D = x_1+\ldots +x_n$ where the points $x_i\in X$ are
not necessarily distinct. More generally a $G$-symmetric product
of $X$ is defined by $SP^n_G(X)= X^n/G$ where $G\subset S_n$ is a
subgroup of the symmetric group on $n$-letters.

If $M$ is a $2$-dimensional manifold, a surface for short, then
$SP^n(M)$ is also a manifold, Section~\ref{sec:everywhere}. In
each theory, relevant examples are essential for illustrating and
understanding general theorems and as a guide for intuition.
Symmetric powers of surfaces provide a list of interesting and
nontrivial examples illustrating many phenomena of manifold
theory.

These manifolds are interesting objects which make surprising
appearance at crossroads of many disciplines of mathematics and
mathematical physics.

If $X=M_g$ is a surface of genus $g$, say a nonsingular algebraic
curve, then $SP^n(M_g)$, interpreted as the space of all effective
divisors of order $n$, serves as the domain of the classical
Jacobi map,\cite{ACGH}, \cite{GrHar},
\[
\mu : SP^n(M_g) \rightarrow {\rm Jac}(M_g) .
\]
In Algebraic Topology, the spaces $SP^{\infty}(X) := {\rm
colim}_{n\in \mathbb{N}}~SP^n(X)$ have, at least for connected
$CW$-complexes, a remarkable homotopical decomposition
\[
SP^{\infty}(X) \simeq \prod_{n=1}^{\infty} K(H_n(X, \mathbb{Z}),
n) ,
\]
due to Dold and Thom, \cite{DoldThom},\cite{Milg_lect}. In
particular the Eilenberg-MacLane space $K(\mathbb{Z}, n)$ has a
natural ``geometric realization'' $SP^\infty(S^n)$. Much more
recent are results which connect special divisor spaces with the
functional spaces of holomorphic maps between surfaces and other
complex manifolds, \cite{Segal}, \cite{C2M2},\cite{Kallel_1998}.
This is a rich theory which can be seen as a part of the evergreen
topological theme of comparing functional spaces with various
particle configuration spaces. It also appears that symmetric
products play more and more important role in mathematical
physics, say in matrix string theory, \cite{Bant2001},
\cite{ByGoOd}, \cite{Dijk}, \cite{DiMoVeVe}, \cite{FujiMatsuo}.

The case of closed, open, orientable or nonorientable surfaces is
as already observed of special interest since their symmetric
products are genuine manifolds. These manifolds were studied in
\cite{Macdonald_topo}, \cite{Mattuck}, \cite{DupontLusztig},
\cite{Arnold_96} and they undoubtedly appeared in many other
papers in different contexts. The signature of $SP^{n}(M_g)$ was
determined by MacDonald \cite{Macdonald_topo}, the signature of
$SP^n(M)$ for more general closed, even dimensional manifolds was
calculated by Hirzebruch \cite{Hirze1968}. Their calculation is
based on the evaluation of the $L$-polynomial and the celebrated
Hirzebruch signature theorem. Zagier in \cite{Zagier1972} used the
Atiyah-Singer $G$-signature theorem and obtained a formula for the
signature of any $G$-symmetric product $SP^n_G(M)$.

The contemporary as well as the classical character of these
objects is particularly well illustrated by Arnold in
\cite{Arnold_96}, where the homeomorphism $SP^n(RP^2)\cong
RP^{2n}$, \cite{DupontLusztig}, is interpreted as a
``quaternionic'' analogue of $SP^n(\mathbb{C}P^1)\cong
\mathbb{C}P^n$ and directly connected with some classical results
of Maxwell about spherical functions.

\medskip
In Section~\ref{sec:everywhere} we give a brief exposition of most
of these results with pointers to relevant references.
Sections~\ref{sec:israel} and \ref{sec:genera} reflect the
research interests of the authors and provide new examples of
applications of symmetric products.

\section{Symmetric products are everywhere!}
\label{sec:everywhere}

\subsection{Examples}
\label{sec:examples}

\begin{defin}
The $n$-th symmetric product or the $n$-th symmetric power of a
space $X$ is $$SP^n(X) = X^n/S_n$$ where $S_n$ is the symmetric
group in $n$ letters.
\end{defin}

Here are first examples of symmetric products of familiar spaces.

 \begin{itemize}
\item[(1)] $SP^{n}(\mathbf{[}0,1])=\Delta ^{n}$ is an $n$-simplex.

\medskip{\small
Since $\mathbf{[}0,1\mathbf{]}$ is totally ordered, $SP^{n}([0,1])=%
\{x_{1}+..+x_{n}\mid 0\leq x_{1}\leq ..\leq x_{n}\leq 1\}=\Delta
^{n}$. The same idea can be used to show that $SP^{n}(\mathbb{R})$
is a closed polyhedral cone in $\mathbb{R}^n$. }

\smallskip
\item[(2)] $SP^{n}(\mathbb{C})=\mathbb{C}^{n}$ and
$SP^{\infty }(\mathbb{C})=\mathbb{C}^{\infty }$.

\medskip{\small
Every element $z_{1}+..+z_{n}\in SP^{n}(\mathbb{C})$ can be
identified with
the monic complex polynomial $p(z)=(z-z_{1})...(z-z_{n})$ with zeros in $%
z_{i}$. Also, $SP^{n}(\mathbb{C})\simeq SP^{n}(\mathbf{\ast
})\cong \ast $.}

\smallskip
\item[(3)] $SP^{n}(S^{1})\simeq S^{1},\,n>0,$ and as a consequence
$SP^{\infty }(S^{1})\simeq S^{1}$.

\medskip{\small
This result follows for example \cite{Milg_lect} from the fact
that the map $\pi :SP^{n}(S^{1})\rightarrow S^{1}$
\[
e^{i\alpha _{1}}+..+e^{i\alpha _{n}}\mapsto e^{i(\alpha
_{1}+..+\alpha _{n})} \]
is a fibration with a contractible fibre.
Alternatively, one can use the homotopy equivalence
$SP^{n}(S^{1})\simeq SP^{n}(\mathbb{C}\setminus \{0 \})$. Then
from $SP^{n}(\mathbb{C})\cong\mathbb{C}^{n}$ we deduce that
$SP^{n}(S^{1})$ is homotopic to the complement of a hyperplane.
\[
SP^{n}(S^{1})\simeq SP^{n}(\mathbb{C}\setminus \{0 \})\simeq
\mathbb{C}^{n}\setminus H\simeq S^{1}. \]
Of course, the Dold-Thom
theorem (Theorem~\ref{thm:DoldThom}) implies $$SP^{\infty
}(S^{1})\simeq \prod_{n=1}^{\infty
}K(H_{n}(S^{1},\mathbb{Z}),n)\simeq
K(H_{1}(S^{1},\mathbb{Z}),1)\simeq K(\mathbb{Z},1)=S^{1}.$$ Note
that $SP^{2}(S^{1})$ is actually the M\"{o}bius band.}

\smallskip
\item[(4)] $SP^{n}(S^{2})=\mathbb{C}P^{n}, \,n>0$ and $SP^{\infty
}(S^{2})=\mathbb{C}P^{\infty}.$

\medskip{\small
First we identify $S^{2}=\mathbb{C}\mathbf{\cup \{\infty \}}$ and
$\mathbb{C} P^{n}=\{p(z)=a_{n}z^{n}+..+a_{0}\mid a_{i}\in
\mathbb{C}\}\diagup (p(z)\sim \lambda p(z),\lambda \neq 0)$. The
map
\[
SP^{n}(S^{2})\ni z_{1}+..+z_{k}+\infty +..+\infty \mapsto
(z+z_{1})..(z+z_{k})\in \mathbb{C}P^{n} \]
is well defined
\begin{eqnarray*}
\lim_{z_{1}\rightarrow \infty }(z+z_{1})..(z+z_{k})
&=&\lim_{z_{1}\rightarrow \infty }
(z^{k}+(z_{1}+..+z_{k})z^{k-1}+..+z_{1}..z_{k}) \\
&=&\lim_{z_{1}\rightarrow \infty }\frac{{1}}{{z}_{1}}
(z^{k}+(z_{1}+..+z_{k})z^{k-1}+..+z_{1}..z_{k})
\\ &=&z^{k}+(z_{1}+..+z_{k-1})z^{k-1}+..+z_{1}..z_{k-1}
\end{eqnarray*}
and easily checked to be a homeomorphism.}

\smallskip
\item[(5)] $SP^{2}(S^{n})=\mathrm{MapCone}(\Sigma ^{n}\mathbb{R}%
P^{n-1}\rightarrow S^{n})$ \cite{Hatcher}$.$
\end{itemize}

\subsection{Maxwell and Arnold}
\label{sec:MaxArn}

The symmetric power of a projective plane is also a projective
space, $SP^n(RP^2)\cong RP^{2n}$, \cite{DupontLusztig}. Vladimir
Arnold in \cite{Arnold_96} observed that this result is a direct
consequence of the theorem on {\em multipole representation of
spherical functions} of James Clerk Maxwell.

Recall that a {\em spherical function} of degree $n$ on a unit
sphere in $\mathbb{R}^3$ is the restriction to the sphere of a
homogeneous harmonic polynomial of degree $n$.

\begin{theo}
\label{thm:Max} The $n$-th derivative of the function
$\frac{1}{r}$ along $n$ constant (translation-invariant) vector
fields $V_1,\ldots V_n$ in $\mathbb{R}^3$ coincides on the sphere
with a spherical function of order $n$. Any nonzero spherical
function $\psi$ of degree $n$ can be obtained by this construction
from some $n$-tuple of nonzero vector fields. These $n$ fields are
uniquely defined by the function $\psi$ (up to multiplication by
nonzero constants and permutation of the $n$ fields).
\end{theo}

In other words any spherical function is the restriction on the
unit sphere of a function of the form
\[ {\cal L}_{V_1}\ldots {\cal L}_{V_n}(\frac{1}{r})\]
where ${\cal L}_X(f)$ is the directional derivative (Lie
derivative) of $f$ in the direction of the vector field $X$.  The
following lemma can be proved by a simple inductive argument based
on the formula
\[ \frac{\partial}{\partial x}\frac{A}{r^b}=
\frac{r^2(\partial A/\partial x)-bAx}{r^b+2}. \]

\begin{lema}
The function ${\cal L}_{V_1}\ldots {\cal L}_{V_n}(\frac{1}{r})$
has the form $P/r^{2n+1}$ where $P$ is a homogeneous polynomial of
degree $n$.
\end{lema}

Actually the polynomial $P$ can be shown to be harmonic in
$\mathbb{R}^3$, i.e. $\Delta(P)= 0$ where $\Delta$ is the Laplace
operator $\Delta = \partial^2/\partial x^2 + \partial^2/\partial
y^2 + \partial^2/\partial z^2$. This essentially follows from the
well known connection between $\Delta$ and the spherical Laplacian
$\tilde\Delta$,
 \begin{equation}
   \label{eqn:laplace}
     \tilde\Delta F = r^2\Delta F -
     \Lambda F, \quad \Lambda := k^2 + k
 \end{equation}
where $F$ is a homogeneous function on $ \mathbb{R}^3\setminus\{
0\}$ of degree $k$, \cite{Arnold_96}. Note that $\tilde\Delta$ is
an operator defined of $k$-homogeneous functions on
$\mathbb{R}^3\setminus\{0\}$, so perhaps a more appropriate
notation would be $\tilde\Delta_k$. From the equation
(\ref{eqn:laplace}) one deduces many important properties of
spherical functions and the associated harmonics in
$\mathbb{R}^3\setminus\{0\}$. For example one immediately observes
that a spherical function, defined as a restriction of a harmonic
function, is an eigen function of the spherical Laplacian.
Conversely, any eigen function $\psi$ of the spherical Laplacian
corresponding to the eigen value $\Lambda = n^2+n$ can be extended
in the ambient space $\mathbb{R}^3\setminus\{0\}$ to a homogeneous
harmonic function in two ways, with the respective degrees $n$ and
$-n-1$ .

Another important consequence is that there is bijective
correspondence between the space $SF_n$ of spherical functions of
order $n$ and the space $HP_n$ of homogeneous, harmonic
polynomials of order $n$. A detailed exposition of these and other
beautiful facts about spherical functions can be found in
\cite{Arnold_96} and \cite{Arnold95}.

Let $S^n(\mathbb{R}^3)$ be the linear space of all homogeneous
polynomials in $\mathbb{R}^3$ of degree $n$. An elementary fact is
that ${\rm dim}(S^k(\mathbb{R}^3))= (n+2)(n+1)/2$. By an
elementary inductive argument one shows that the linear map
$\delta_n : S^n(\mathbb{R}^3)\rightarrow  S^{n-2}(\mathbb{R}^3)$
defined as the restriction of the Laplacian $\Delta$, is an
epimorphism. It follows that the dimension of the space $HP_n =
{\rm Ker}(\delta_n)$ is $(n+2)(n+1)/2 - n(n-1)/2 = 2n+1$ so the
dimension of $SF_n$ is also $2n+1$. Finally by
Theorem~\ref{thm:Max}, one observes that the symmetric power of
$RP^2$ is homeomorphic to the projective space associated to
$SF_n$, hence
\[ SP^n(RP^2) \cong RP^{2n}.\]

\subsection{Abel and Jacobi}

Suppose that $M=M_g$ is a compact Riemann surface, i.e. a compact
$2$-manifold (of genus $g$) with a complex structure.
Alternatively, $M$ can be introduced as a nonsingular algebraic
curve. Symmetric powers $SP^n(M)$ of a curve are ubiquitous in
Algebraic Geometry, \cite{Forster} \cite{GrHar} \cite{Har}. We
cannot possibly do justice to most of these developments in this
article. Keeping in mind our focus on the {\em topological
manifold} $SP^n(M_g)$, we start with the following result
\cite{Mattuck} which is a generalization of the fact
$SP^n(\mathbb{C}P^1)\cong \mathbb{C}P^n$.

\begin{theo}
 \label{thm:fibr}
 Suppose that $n> 2g-2$. Then there is a fibre bundle
\begin{equation}\label{eqn:fibration}
  \mathbb{C}P^{n-g} \rightarrow SP^n(M_g) \rightarrow T^{2g}
\end{equation}
where the fibre $\mathbb{C}P^{n-g}$ is a complex projective space
of dimension $n-g$ and $T^{2g}$ is a $2g$-dimensional torus.
\end{theo}

We start an outline of the proof of this theorem with a brief
exposition of the Abel--Jacobi map.

A standard fact \cite{Forster} is that the complex vector space
$\Omega(M_g)$ of {\em holomorphic} differential $1$-forms is
$g$-dimensional. Let $\omega_1, \ldots ,\omega_g$ be a basis of
this space. Define the {\em subgroup of periods} ${\rm Per} = {\rm
Per}(\omega_1,\ldots, \omega_g)$ in $\mathbb{C}^g$ by the
requirement that $v = (v_1,\ldots, v_g)\in {\rm Per}$ if and only
if for some $\alpha\in\pi_1(M_g)$
\[ \mbox{ \rm for all } i=1,\ldots, g \quad v_i = \int_\alpha\,\omega_i \, .\]
Then ${\rm Per}$ is a discrete subgroup in $\mathbb{C}^g$ of
maximal rank, hence ${\rm Jac}(M_g):= \mathbb{C}^g/{\rm Per}$ is a
$(2g)$-dimensional torus called the {\em Jacobian} of the surface
$M_g$.

Suppose that $b\in M_g$ is a base point. Given $x\in M_g$ and a
path $\beta$ connecting points $b$ and $x$, let $u = (u_1,\ldots
,u_g)$ be a vector in $\mathbb{C}^g$ defined by
\[ u_i := \int_\beta\, \omega_i\, .\]
The vector $u\in \mathbb{C}^g$ depends on $\beta$, however its
image in ${\rm Jac}(M_g)$ depends only on the point $x$. This way
arises the celebrated {\em Abel--Jacobi} map
 \begin{equation}
 \label{eqn:AbJac} \mu : M_g\rightarrow {\rm Jac}(M_g).
 \end{equation}
If $g=1$, i.e. in the case of an {\em elliptic curve}, the map
$\mu$ is an isomorphism. This is a famous turning point in
mathematics, when the study of meromorphic functions on the curve
was reduced to the study of meromorphic functions on a torus, i.e.
the meromorphic functions in the complex plane $\mathbb{C}$ with
$2$ periods.

In the case of a general curve $M_g$, the map $\mu$ is far from
being an isomorphism. Since ${\rm Jac}(M_g)$ is an abelian group,
the Abel--Jacobi map $\mu$ can be extended to a symmetric power
$SP^n(M_g)$ by a formula
\[ \mu_n(D) := \mu(x_1) + \ldots + \mu(x_n)\]
where $D = x_1 +\ldots + x_n \in SP^n(M_g)$. If $n=g$ then $\mu_g
: SP^g(M_g)\rightarrow T^{2g}$ is a ``correct'' replacement for
the map $\mu$.

\medskip
An alternative description of the map $\mu_n$ is following. Let
${\rm Pic(M_g) := {\rm Div}_0/{\rm Div}_H}$ be the {\em Picard}
group of $M_g$ where ${\rm Div}_0$ is the group of all divisors of
degree $0$ and ${\rm Div}_H$ is the group of principal divisors,
i.e. the divisors of the form $D = (f)$ for some meromorphic
function $f$. There is a map
\[ \Phi : {\rm Div}_0\rightarrow {\rm Jac}(M_g)\]
defined by $\Phi(D) = v = (v_1,\ldots, v_g)$ where
\[ v_i := \int_c\,\omega_i \quad \mbox{ {\rm for each } } i \]
and $c$ is a $1$-chain, i.e. a system of paths connecting points
in $D$, such that $\partial(c) = D$. Abel's theorem \cite{GrHar}
claims that the kernel of $\Phi$ is precisely the group ${\rm
Div}_H$ of principal divisors, hence the induced map
\begin{equation}
\label{eqn:map} \phi : {\rm Pic}(M_g) \rightarrow {\rm Jac}(M_g)
\end{equation}
is a monomorphism. It turns out\footnote{Jacobi inversion problem,
\cite{GrHar}.} that the map $\phi$ is an isomorphism. As a
consequence, the Abel--Jacobi map $\mu=\mu_1$ (\ref{eqn:AbJac}),
more precisely its higher dimensional extension $\mu_n$, has a
twin map
\begin{equation}
\label{eqn:twin} \nu_n : SP^n(M_g)\rightarrow {\rm Pic}(M_g)
\end{equation}
defined by $\nu_n(D) = D - nb$ ($b$ is the base point in $M_g$),
such that $\phi\circ \nu_n = \mu_n$.

\medskip
As a consequence one can approach the proof of
Theorem~\ref{thm:fibr} via the map $\nu_n$. In order to determine
the inverse image $\nu_n^{-1}([D])$ for a given $[D]\in {\rm
Pic}(M_g)$, one has to solve the Riemann--Roch problem i.e. to
determine the dimension of the space of effective divisors with
prescribed information about their zeros and poles. By the
Riemann--Roch theorem, if $n>2g-2$ then the dimension of the space
of meromorphic functions $f$ such that $E = D + (f)$ is an
effective divisor of order $n$ is precisely $n-g+1$. The
representation of the divisor $E$ in the form $D+(f)$ is not
unique, namely if $D+(f) = D+(g)$ then $(f)=(g)$ and $g = cf$ for
some nonzero constant $c\in \mathbb{C}$. It follows that
$\nu_n^{-1}(D)$ is a complex projective space of dimension $n-g$
which finally explains the appearance of the fibre
$\mathbb{C}P^{n-g}$ in the fibration (\ref{eqn:fibration}).

\subsection{The Poincar\' e polynomial of a symmetric product}
\label{sec:Mac}

In this section we compute the Betti numbers of general symmetric
products. These results were originally obtained by
I.G.~MacDonald, \cite{Macdonald_62}.

\medskip
 Let $V=\oplus_{d\geq 0}V_d$ be a graded
finite dimensional vector space. The associated Poincar\' e
polynomial is defined by $ P_t(V)=\sum_{d\geq 0}t^d {\rm dim}V_d
.$ It is easily shown that
\[ P_t(V\oplus W)=P_t(V)+P_t(W),\]
\[ P_t(V\otimes W)=P_t(V)P_t(W).\]
The symmetric algebra over the vector space $V$ is defined by
$$S^*(V)=T^*(V)/v\otimes w-(-1)^{degv\cdot degw}w\otimes v .$$ It
is naturally bigraded by $${\rm bdeg} \langle v_1,\ldots,v_n
\rangle =(\sum_{j=1}^n{\rm deg}(v_j),n).$$ We introduce the formal
variable $q$ by $S_q^*(V)=\sum_{n\geq 0}q^nS^n(V),$ where
$S^n(V)=\{x\in S^*(V) \mid {\rm bdeg}(x)=(\cdot,n)\}$ is the
$n-$th symmetric power of $V$. Since $ S^n(V\oplus
W)=\oplus_{p+q=n}S^p(V)\otimes S^q(W),$ we have
\[ S_q^*(V\oplus W)=S_q^*(V)S_q^*(W),\]
\[ P_t(S_q^*(V\oplus W))=P_t(S_q^*(V))P_t(S_q^*(W)).\]
If $L= L_d$ is a $1$-dimensional, graded vector space generated by
a vector of degree $d$ then,
\[
 S_q^*(L_d)=\left\{
\begin{array}{ll}1+qL+q^2L^{\otimes 2}+\cdots & \mbox{, d even}\\
1+qL & \mbox{, d odd}
\end{array}\right. , \mbox{ { \rm and } } \]

\[ P_t(S^*_q(L_d))=\left\{
\begin{array}{ll} \frac{1}{1-qt^d} & \mbox{ , d even}\\ 1+qt^d &
\mbox{ , d odd}  \end{array} \right. .
 \] Since any graded vector space
$V\oplus_{d\geq 0}V_d$ is decomposable into a sum of
$1$-dimensional, graded vector spaces, it follows that,
\[ P_t(S_q^*(V))=\frac{\prod_{d\,
odd}(1+qt^d)^{dimV_d}}{\prod_{d\, even}(1-qt^d)^{dimV_d}}.\] Now,
if $V=H_*(X, \mathbb{Q})$ is the homology of a CW-complex $X$ and
$\beta_d$'s are its Betti numbers, then \[ H_*(SP^n(X),
\mathbb{Q})=(H_*(X, \mathbb{Q})^{\otimes n})^{S_n}=S^n(H_*(X,
\mathbb{Q})=S^n(V)\] and we get the MacDonald result
\cite{Macdonald_62}
\[ \sum_{n\geq 0}q^nP_t(SP^n(X))=\frac{\prod_{d\,
odd}(1+qt^d)^{\beta_d}}{\prod_{d\, even}(1-qt^d)^{\beta_d}}.\]

In particular for $t=-1$ we get the generating function for Euler
characteristics of symmetric powers of the space $X$
\[ \sum_{n\geq 0}q^n\chi(SP^n(X))=(1-q)^{-\chi(X)}.\]

\subsection{Dold-Thom theorems}

Symmetric powers of spaces have been used in homotopy theory for
the last fifty years. They for example appear in the study of
iterated loop spaces (the symmetric products appear e.g. as
fragments of model spaces for important spaces such as $\Omega
^{n}\Sigma ^{n}X$). Here we review some of the central results.

\medskip
Let $(X,\ast )$ be a space with the base point $\ast \in X$.
Assuming that $SP^{0}(X)=\{\ast \}$, for each $n\geq 0$ we define
a natural inclusion
\[ SP^{n}(X)\hookrightarrow SP^{n+1}(X), \quad
x_{1}+...+x_{n}\mapsto x_{1}+...+x_{n}+\ast .
\]

The colimit of the direct system of these inclusions is the so
called infinite symmetric product $SP^{\infty }(X)$. The ``we can
always add'' map $\mu :SP^{n}(X)\times SP^{m}(X)\rightarrow
SP^{n+m}(X)$
\[
(x_{1}+...+x_{n},y_{1}+...+y_{m})\mapsto
x_{1}+...+x_{n}+y_{1}+...+y_{m}
\]
induces associative multiplication on $SP^{\infty }(X)$ with the
neutral element $\ast \in SP^{0}(X)$. Moreover, it can be proved
that $SP^{\infty }(X)$ is the free commutative topological monoid
generated by $X$ with $\ast $ as the unit element. It may be
natural to ask what can be said about $A^{\infty }(X),$ the free
commutative topological group generated by $X$ with $\ast $ as
neutral and the topology of the quotient

\[
A^{\infty }(X)=\coprod_{n,m\geq 1} SP^{n}(X)\times SP^{m}(P)/ \sim
\]
where the equivalence relation $\sim$ is defined by
\[
(x_{1}+..+x_{n},y_{1}+..+y_{m})\sim
(x_{1}+..+\hat{x}_{i}+..+x_{n},y_{1}+..+ \hat{y}_{j}+..+y_{m})
\]
if and only if $x_{i}=y_{j}$.

 A nonzero element in $A^{\infty }(X)$ can be
formally written as a difference
$(x_{1}+..+x_{n})-(y_{1}+..+y_{m})$ where elements $x_{i}$,
$y_{j}$ are all different from the base point $\ast $.

\begin{theo}{\rm( Dold-Thom)}
\label{thm:DoldThom} If $(X,\ast )$ is a connected $CW$-complex,
then
\[
SP^{\infty }(X)\simeq \prod_{n=1}^{\infty
}K(H_{n}(X,\mathbb{Z}),n)
\]
where $K(G,n)$ is an Eilenberg-MacLane space, i. e. a $CW$-complex
with the property that $\pi _{n}(K(G,n))=G$ and $\pi
_{i}(K(G,n))=0$ for each $i\neq n$.
\end{theo}

{\bf Examples:}

\begin{itemize}

\item[(1)] \qquad $SP^{\infty }(S^{n})\cong K(\mathbb{Z},n)$.

\item[(2)] \qquad $SP^{\infty }(\mathbb{C}P^{n})=
\prod_{k=1}^n K(\mathbb{Z},2k)= \prod_{k=1}^n SP^{\infty
}(S^{2k})$.

\item[(3)] \qquad $SP^{\infty }(S^{n}\cup
_{k}e^{n+1})=K(\mathbb{Z}/k,n)$.

\end{itemize}
There are different proofs of the theorem of Dold and Thom, see
\cite{DoldThom},\cite{Milg_lect}. One possibility is to establish
first the following relative of Theorem~\ref{thm:DoldThom}.

\begin{theo}
If $(X,\ast )$ is a connected $CW$-complex, then
\[
A^{\infty }(X)\simeq \prod_{n=1}^{\infty
}K(\tilde{H}_{n}(X,\mathbb{Z}),n).
\]
\end{theo}

The proof of this theorem given in \cite{Milg_lect} is based on
the following facts:

\begin{itemize}

\item[(i)] $X\simeq Y$ $\Longrightarrow A^{\infty }(X)\simeq
A^{\infty }(Y)$ \\(a homotopy $H:X\times I\rightarrow Y$ yields a
homotopy $A^{\infty }(H):A^{\infty }(X)\times I\rightarrow
A^{\infty }(Y)$).

\item[(ii)] $A^{\infty }(S^{0})\cong \mathbb{Z}$.

\item[(iii)] $A^{\infty }(S^{n})=K(\mathbb{Z},n)$,
(a cofibration sequence $S^{n}\hookrightarrow D^{n+1}\rightarrow
S^{n+1}$ produces a
fibration $A^{\infty }(D^{n+1})\rightarrow A^{\infty }(S^{n+1})$ with $%
A^{\infty }(S^{n})$ as a fibre, so an induction on the dimension
can be applied).

\item[(iv)] $X\mapsto \pi _{i}(A^{\infty }(X))$ induces a reduced
homology theory with integral coefficients (it satisfies
Eilenberg-Steenrod axioms: (i), $A^{\infty }(\ast )\simeq \ast $,
(iii)).

\item[(v)] The uniqueness of the homology theory satisfying the Eilenberg-Steenrod
axioms implies that $\pi_{i}(A^{\infty }(X))\cong
\tilde{H}_{i}(X,\mathbb{Z})$.

\end{itemize}

\begin{theo}{\rm (Dold-Thom)}
\label{thm:druga}
If $(X,\ast )$ is a connected $CW$-complex then
the inclusion $SP^{\infty
}(X)\hookrightarrow A^{\infty }(X)$%
\[
x_{1}+...+x_{n}\mapsto x_{1}+...+x_{n}
\]
is a homotopy equivalence.
\end{theo}

By using the comparison theorem for spectral sequences, one can
prove Theorem~\ref{thm:druga} for spheres. Then adding cell after
cell, with the use of the 5-lemma, the theorem is established for
every connected $CW$-complex.

As a consequence one concludes that both $X\mapsto SP^{n}(X)$ and
$X\mapsto SP^{\infty }(X)$ are {\em homotopy contravariant
functors}.

\begin{cor}
$SP^{\infty }(X\vee Y)\simeq SP^{\infty }(X)\times SP^{\infty
}(Y)$.
\end{cor}

This is a direct consequence of the Dold-Thom theorem
(Theorem~\ref{thm:DoldThom}) and the following facts
$\tilde{H}_{n}(X\vee Y,\mathbb{Z})\cong \tilde{H}_{n}
(X,\mathbb{Z})\oplus \tilde{H}_{n}(Y,\mathbb{Z}),$ $K(G\times
H,n)=K(G,n)\times K(H,n)$.

\subsection{Steenrod and Milgram;\\ homology of symmetric products}

The inclusion map $i:SP^{n}(X)\hookrightarrow SP^{n+1}(X)$ is very
useful in computations with symmetric products. For example this
map induces a long exact sequence in homology (with any
coefficient group)
\begin{equation}
\label{eqn:long} ..\rightarrow H_{\ast
}(SP^{n}(X))\overset{i_{\ast }}{\rightarrow }H_{\ast
}(SP^{n+1}(X))\rightarrow H_{\ast
}(SP^{n+1}(X),SP^{n}(X))\rightarrow ..
\end{equation}
Consequently it is essential to understand the associated map
$i_{\ast }$. A solution of this problem was announced by Norman
Steenrod in \cite{Steenrod} but this proof has never been
published. Albrecht Dold proved in  \cite{Dold} the following
theorem.

\begin{theo}
\label{thm:Steenrod} Let $(X,\ast )$ be a connected $CW$-complex,
$G$ an arbitrary coefficient group
and $i_{n}:SP^{n}(X)\hookrightarrow SP^{n+1}(X)$, $j_{n}:SP^{n}(X)%
\hookrightarrow SP^{\infty }(X)$ the natural inclusions. Then
$(i_{n})_{\ast }$ is an inclusion onto a direct summand, i. e.
there is a splitting exact sequence
\[
0\rightarrow H_{n}(SP^{n}(X),G)\overset{(i_{n})_{\ast }}{\rightarrow }%
H_{n}(SP^{n+1}(X),G) .
\]
\end{theo}

In light of Theorem~\ref{thm:Steenrod}, the long exact sequence
(\ref{eqn:long}) becomes
\[
..\overset{0}{\rightarrow }H_{\ast }(SP^{n}(X))\overset{i_{\ast }}{%
\rightarrow }H_{\ast }(SP^{n+1}(X))\overset{onto}{\rightarrow
}H_{\ast }(SP^{n+1}(X),SP^{n}(X))\overset{0}{\rightarrow }..
\]
and implies
\begin{eqnarray*}
H_{\ast }(SP^{n}(X),G) &=&\bigoplus_{i=1}^n H_{\ast
}(SP^{i}(X),SP^{i-1}(X),G)\text{, } \\ H_{\ast }(SP^{\infty
}(X),G) &=&\bigoplus_{i=1}^\infty H_{\ast
}(SP^{i}(X),SP^{i-1}(X),G).
\end{eqnarray*}
 Here we used the commutativity of the following diagram
\[
\begin{array}{c}
H_{n}(SP^{n}(X),G)\overset{(i_{n})_{\ast }}{\longrightarrow }%
H_{n}(SP^{n+1}(X),G) \\
^{(j_{n})_{\ast }}{ \searrow }\text{ \ \ \ \ }{ \swarrow }%
^{(j_{n+1})_{\ast }} \\ H_{n}(SP^{\infty }(X),G)
\end{array}
\]
This suggests that there should exist a natural filtration $$
H_{1}\subseteq H_{2}\subseteq ..\subseteq H_{n}\subseteq ..H_{\ast
}(SP^{\infty }(X),G)$$ of $H_{\ast }(SP^{\infty }(X),G)$ such that
$H_{n}\diagup H_{n-1}\cong H_{\ast }(SP^{n}(X),SP^{n-1}(X),G)$.
Since the homology  is compactly supported, for each $\alpha \in
H_{m}(SP^{\infty }(X),G)$  we define
\[
n_{\alpha }=\min \{r\mid (\exists \alpha _{r}\in
H_{m}(SP^{r}(X),G))\,i_{r}(\alpha _{r})=\alpha \}.
\]
Thus, there is a filtration
\[
H_{n}=\{\alpha \in H_{\ast }(SP^{\infty }(X),G)\mid n_{\alpha
}\leq n\}
\]
and it is obvious that $H_{n}\subseteq H_{n+1}$\ and
\[
H_{n}\diagup H_{n-1}=\{\alpha \in H_{\ast }(SP^{\infty }(X),G)\mid
n_{\alpha }=n\}\cong H_{\ast }(SP^{n}(X),SP^{n-1}(X)).
\]
Also, every filtration element $H_{n}$ is additionally filtered
with the groups
\[
F_{n,m}=H_{m}(SP^{\infty }(X),G)\cap H_{n}.
\]

Now the ''we can always add'' map $\mu :SP^{n}(X)\times
SP^{m}(X)\rightarrow SP^{n+m}(X)$ induces a Pontriagin product on
filtration elements
\[
\mu :F_{n,m}\otimes F_{i,j}\rightarrow F_{n+i,m+j}\text{.}
\]
Hence, for ``untwisted'' coefficients, say for a field
$\mathbb{K}$,
\begin{eqnarray*}
H_{\ast }(SP^{\infty }(X),\mathbb{K)} &=& \bigoplus_{i=1}^\infty
H_{\ast }(SP^{i}(X),SP^{i-1}(X),\mathbb{K})
\\
&=&\bigoplus_{i=1}^\infty \bigoplus_{j=1}^\infty
H_{j}(SP^{i}(X),SP^{i-1}(X),\mathbb{K})= \bigoplus_{i=1}^\infty
\bigoplus_{j=1}^\infty F_{i,j}\diagup F_{i-1,j}
\end{eqnarray*}
is a bigraded, commutative associative algebra with a neutral
element.

\bigskip

James Milgram in \cite{Milg} gave an another idea for calculating
the homology of symmetric product $SP^{\infty }(X)$. The first
step is a homology decomposition of the space $X$ in a wedge
$\bigvee_{i\in I}M_{i}$ of Moore spaces. So instead of $H_{\ast
}(SP^{\infty }(X))$ we calculate $H_{\ast }(SP^{\infty
}(\bigvee_{i\in I}M_{i}))$. Knowing that $SP^{\infty }(X\vee
Y)\simeq SP^{\infty }(X)\times SP^{\infty }(Y)$ it can be proved
that there is bigraded algebra isomorphism
\[
H_{\ast }(SP^{\infty }(X\vee Y),\mathbb{K})\cong H_{\ast }(SP^{\infty }(X),%
\mathbb{K})\otimes H_{\ast }(SP^{\infty }(Y),\mathbb{K}).
\]
Hence,
\[
H_{\ast }(SP^{\infty }(X),\mathbb{K})\cong H_{\ast }(SP^{\infty
}(\bigvee_{i\in I}M_{i}),\mathbb{K})\cong \bigotimes_{i\in
I}H_{\ast }(SP^{\infty }(M_{i}),\mathbb{K}).
\]
So it remains to determine $H_{\ast }(SP^{\infty }(M))$ for a
Moore space $M$. Recall that $M$ is a Moore space of type $(G,n)$
if $M$ is $CW$-complex with one $0$-cell and other cells only in
dimensions $n$ and $n+1$, such that $H_{n}(M)=G $ and
$\tilde{H}_{i}(M)=0$ for $i\neq n$. According to Dold-Thom theorem
\[
H_{\ast }(SP^{\infty }(X),\mathbb{K})\cong H_{\ast
}(K(G,n),\mathbb{K}).
\]
Finally, $H_{\ast }(K(G,n),\mathbb{K})$ is determined from the
spectral sequence of the fibre space
\[
\begin{array}{cc}
K(H,n)\longrightarrow & K(G,n) \\ &
\begin{array}{c}
{ \downarrow } \\ K(K,n)
\end{array}
\end{array}
\]
where $0\rightarrow H\rightarrow G\rightarrow K\rightarrow 0$ is
exact sequence of abelian groups.

For example if $X$ is compact Riemann surface of genus $g$ there
is a homology isomorphism $H_{\ast }(X,\mathbb{Z})\cong H_{\ast
}(\bigvee_{1}^{2g}S^{1}\vee S^{2})$ and so
\[
H_{\ast }(SP^{\infty }(X),\mathbb{Z})\cong \Lambda (e_{1})\otimes
..\otimes \Lambda (e_{2g})\otimes \Gamma \lbrack f]
\]
where $\deg e_{i}=1$, $e_{i}\in F_{1,1}$, $\deg f=2$, $f\in F_{1,2}$ and $%
\Gamma \lbrack f]$ is the divided power algebra.

\medskip\noindent
{\bf Question:} What can one say about the group $H_{\ast
}(SP^{\infty }(X),\mathbb{Z})$ in the case of a simply connected
$4$-manifold $X$?

\subsection{Divisor spaces}

Graeme Segal studied in \cite{Segal} the topology of rational
functions of the form $f = p(z)/q(z)$ where $p$ and $q$ are monic
polynomials of degree $n$ which do not have a common root. The
space $F_n$ of such polynomials can be identified with a subspace
of the space $M_n$ of all self maps of the Riemann sphere $S^2 =
\mathbb{C}\cup\{\infty\}$ of degree $n$ which take $\infty$ to
$1$. The first Seagal's result is that the inclusion
$F_n\rightarrow M_n$ is a homotopy equivalence up to dimension
$n$.

This and other results of Seagal served as a motivation for
J.~Milgram \cite{Milg_lect} and other mathematicians to study the
the so called {\em divisor spaces} or {\em particle space} as
approximations for important, infinite dimensional functional
spaces ${\rm Map}(U,V)$ and their subspaces. For a $2$-dimensional
complex manifold $X$ one defines the divisor space
\begin{equation*}
\mathrm{Div}_{k}(X)=\{(\Sigma _{i=1}^{k}x_{i},\Sigma
_{i=1}^{k}y_{i})\in SP^{k}(X)\times SP^{k}(X)\mid
\{x_{1},..,x_{k}\}\cap \{y_{1},..,y_{k}\}=\emptyset \}.
\end{equation*}

Kallel \cite{Kallel_1998} introduces even more general spaces
\begin{equation*}
\mathrm{Div}_{k_{1},..k_{n}}(X)=\{(\Sigma
_{i=1}^{k_{j}}x_{ij})_{j=1}^{n}\in
\prod_{j=1}^{n}SP^{k_{j}}(X)\mid
\bigcap_{j=1}^{n}\{x_{1j},..,x_{k_{j}j}\}=\emptyset \}.
\end{equation*}

By building appropriate model and application of spectral sequence
he proved the following result
\begin{equation*}
H_{\ast }(\mathrm{Div}_{k_{1},..k_{n}}(X\backslash
\{x\})\mathbb{K},)\cong \mathrm{Tor}_{2nk-\ast ,k}^{H_{\ast
}(SP^{\infty }(X))}(\mathbb{K},H_{\ast }(SP^{\infty
}(X);\mathbb{K})^{\otimes n})
\end{equation*}
for the Riemann surface $X$ of genus $g$ and a coefficient field $\mathbb{K%
}$.

\subsection{Dupont and Lusztig}

Suppose that $X$ is a compact, closed {\em unorientable}
$2$-manifold such that \begin{equation}\label{eqn:dim} {\rm
dim}(H_1(X, \mathbb{Q})) = g\, .\end{equation}
The following
theorem was proved by J.L.~Dupont and G.~Lusztig in
\cite{DupontLusztig},

\begin{theo}
\label{thm:DupLus} For $n\geq g$, the symmetric product $SP^n(X)$
is diffeomorphic to a $(2n-g)$-dimensional real projective bundle
over the $g$-dimensional torus $T^g$.
\end{theo}

The proof of Dupont and Lusztig follows the idea of the proof of
Theorem~\ref{thm:fibr}. The starting point is an observation that
for a given unorientable $X$, satisfying the condition
(\ref{eqn:dim}), there exists a Riemann surface $Y$ together with
a fixed point free, antiholomorphic involution $T : Y\rightarrow
Y$ such that $X\cong Y/(\mathbb{Z}/2)$ with the
$\mathbb{Z}/2$-action determined by $T$. Then the genus of $Y$ is
$g$. Let ${\cal J}_{2n}$ be the set of all isomorphism classes of
holomorphic line bundles on $Y$ with the Chern class equal to
$2n$. Then ${\cal J}_{2n}$ is a free homogeneous space of the
complex torus (identified as ${\cal J}_0$), of complex dimension
$g$.

Let $L_y$ be the holomorphic line bundle associated to the divisor
$y$ for some $y\in Y$. Following \cite{Macdonald_topo}, the map
\[ (y_1,\ldots, y_{2n}) \mapsto L_{y_1}\otimes\ldots\otimes
L_{y_{2n}}\] defines for $n\geq g$ a holomorphic projective bundle
over ${\cal J}_{2n}$,

\begin{equation}\label{eqn:bun}
\begin{CD}
 \mathbb{C}P^{2n-g}\longrightarrow SP^{2n}(Y) @> \Phi>>
{\cal J}_{2n}
\end{CD}
\end{equation}
The involution $T$ acts antiholomorphically on $SP^{2n}(Y)$ by the
formula $$T(y_1,\ldots, y_{2n}) = (Ty_1,\ldots, Ty_{2n})$$ and on
${\cal J}_{2n}$ by the formula $T(L) = T^\ast(\bar{L})$. It is
easily verified that the map $\Phi$ in (\ref{eqn:bun}) is
equivariant with respect to this action. Hence it takes one fibre
of this map into another in an antiholomorphic way preserving
their projective structures.

One deduces from here that the fixed point set $A$ of $T:
SP^{2n}(Y)\rightarrow SP^{2n}(Y)$ is a real projective bundle over
a union of components, denoted by ${\cal J}_n'$, of the fixed
point set $B$ of $T: {\cal J}_{2n}\rightarrow {\cal J}_{2n}$.
Using the fact that $T: Y\rightarrow Y$ is fix-point free, one
sees that $A$ can be identified with $SP^n(Y/T) = SP^n(X)$, hence
in particular it is connected. The space $B$ is a free homogeneous
space of the subgroup of ${\cal J}_0$ fixed by $T$, hence a union
of $g$-dimensional tori. It follows that ${\cal J}_n'$ is a real
$g$-dimensional torus and the theorem is proved.

\subsection{New invariants of $3$--manifolds}

In this section we briefly outline the role of symmetric powers of
Riemann surfaces in a recent progress \cite{OzsSza} in
constructing invariants of $3$-manifolds via Floer homology. Floer
\cite{Floer} originally defined his groups for a symplectic
manifold $(M,\omega)$ and a pair $\Sigma_1$ and $\Sigma_2$ of its
Lagrangian submanifolds. P.~Ozsv\' ath and Z.~Szab\' o show in
\cite{OzsSza} how a similar theory, producing new invariants of
$3$-manifolds, can be developed with the symmetric power
$SP^g(M_g)$ in the role of the symplectic manifold $M$.

Recall that a Heegaard splitting of a $3$-manifold $U$ is a
decomposition of the form $U = U_0\cup_\Sigma U_1$, where $U_0$
and $U_1$ are two handlebodies glued together by an orientation
preserving diffeomorphism $\phi : \partial(U_0)\rightarrow
\partial(U_1)$ of their boundaries. The common boundary $\Sigma$
is assumed to be a Riemann surface $M_g$ of genus $g$. The isotopy
class of the diffeomorphism $\phi$, and the associated Heegaard
splitting, are determined by two collections
$\{\alpha_1,\ldots,\alpha_g\}$ and $\{\beta_1,\ldots,\beta_g\}$ of
simple, closed curves in $\Sigma=M_g$. The symmetric power
$SP^g(\Sigma)$ is a complex manifold, and if
$\langle\cdot,\cdot\rangle$ is an associated Hermitian metric and
$J$ its complex structure, then the $2$-form $\omega$ defined by
$\omega(X,Y):= \langle X, JY\rangle$, (if closed) turns
$SP^g(\Sigma)$ into a symplectic manifold. Actually Ozsv\' ath and
Szab\' o show how Floer homology groups, reflecting the properties
of the input $3$-manifold $U$, can be defined with the torii
$T_1=\alpha_1\times\ldots\times\alpha_g$ and
$T_2=\beta_1\times\ldots\times\beta_g$ in the role of Lagrangian
submanifolds $\Sigma_1$ and $\Sigma_2$. These are in general only
totally real submanifolds of $SP^g(\Sigma)$ and in order to define
Floer homology groups they need an additional hypothesis that $U$
carries a structure of a ${\rm Spin}^c$ manifold.

This is just a beginning of a beautiful  and interesting theory
and the reader is referred to \cite{OzsSza} and subsequent
publications for details.

\section{Symmetric powers of open surfaces}
\label{sec:israel}

\subsection{Signature of symmetric products of punctured surfaces}

\begin{defin} Given a Riemann surface of genus $g$, the associated
{\em open} or {\em punctured} surface $M_{g,k}$ is defined by
\[M_{g,k} = M_g\setminus\{\alpha_1,\ldots,\alpha_k\}\]
where $\{\alpha_1,\ldots,\alpha_k\}$ is a collection of $k$
distinct points in $M_g$.
\end{defin}

One of the main results of \cite{BGZ01} is the proof of the
existence of punctured surfaces $M = M_{g,k}$ and $N = M_{g',k'}$
such that the associated symmetric products $SP^{2n}(M)$ and
$SP^{2n}(N)$ are not homeomorphic although $M$ and $N$ have the
same homotopy type. Actually it was shown that this is the case
for the punctured surfaces $M_{g,k}$ and $M_{g',k'}$ which satisfy
the conditions:
\begin{itemize}
\item $2g + k = 2g' + k'$,
\item $g\neq g'$ and ${\rm max}\{g,g'\}\geq n$.
\end{itemize}
The key ingredient in the proof that the associated symmetric
products $SP^{2n}(M_{g,k})$ and $SP^{2n}(M_{g',k'})$ are not
homeomorphic is the computation of the signature of these (open!)
manifolds.

\begin{theo}{\rm (\cite{BGZ01})}
\label{thm:sign-C}
\begin{equation}
\label{sign-3} {\rm Sign}(SP^{2n}(M_{g,k})) =
(-1)^n{\binom{g}{n}}\, .
\end{equation}
\end{theo}

\subsection{A connection with $(m+k,m)$-groups}

An $(m+k,m)$-grupoid or a vector valued grupoid $(G,f)$ is simply
a map $f : G^{m+k}\rightarrow G^m$. The analogs of commutativity,
associativity and other algebraic laws can be formulated for these
objects and the corresponding algebraic structures are called
$(m+k,m)$-semigroups, $(m+k,m)$-groups etc.

The theory of vector valued algebraic structures was developed in
the eighties by G.~\v Cupona, D.~Dimovski, K.~Tren\v cevski and
their collaborators, see \cite{TreDim}, \cite{TreDim1} and the
references in \cite{TreDim}. Perhaps a motivation for the study of
these objects, aside from the intrinsic algebraic interest, can be
found in a growing interest in vector valued structures following
the development of the {\em theory of operads}, \cite{Ope}.

\medskip
Our point of departure is an observation\footnote{We are indebted
to Prof. Kostadin Tren\v cevski for the information that the
symmetric products of surfaces are relevant for the theory of
$(m+k,m)$-groups.} that if $(M,f)$ is a topological, commutative
$(m+k,m)$-group, then the symmetric product $SP^m(M)$ admits the
structure of a commutative Lie group.

\medskip
Here are  relevant excerpts from \cite{TreDim}:

\begin{itemize}

\item
\noindent {\bf Theorem 6.1.} ({\em {\rm \cite{TreDim1}} Theorem
3.5 $\&$ Prop. 3.2 Chap. III) If $(M,f)$ is locally euclidean
topological, commutative $(m+k,m)$-group for $m\ge 2$, then $dim
(M)=2$, $M$ is oriented manifold not homeomorphic to the sphere
$S^{2}$ and } $$SP^m(M) \cong  {\bf R}^{u} \times (S^{1})^{v}.$$

\item {\bf Conjecture.} {\em Each connected, locally euclidean
topological, commutative $(m+k,m)$-group is isomorphic to an
affine $com(m+k,m)$-group.}

\end{itemize}

It is clear that the results from Section~\ref{sec:israel} are
relevant for this conjecture. In other words the signature
computation (Theorem~\ref{thm:sign-C}) rule out many open surfaces
$M$ as possible ground spaces for a  structure of a locally
euclidean $(m+k,m)$-group. A more detailed analysis will be
published in a subsequent publication.

\section{Genera of symmetric powers}
\label{sec:genera}

Let $X$ be a closed, oriented manifold with orientation preserving
action of a finite group $G$. Rational cohomology of the orbit
space $X/G$ is naturally isomorphic to $G-$invariant part of
rational cohomology of $X$. The equivariant Euler characteristic
is defined as \[\chi(g,X)=\sum_{j\geq0}(-1)^j{\rm
tr}g^*|_{H^j(X)}.\] It is easy fact from representation theory
that
\[\chi(X/G)=\frac{1}{|G|}\sum_{g\in G}\chi(g,X).\] The same
formula holds for signatures of manifolds \cite{Zagier1972} \[{\rm
sign}(X/G)=\frac{1}{|G|}\sum_{g\in G}{\rm sign}(g,X).\]

Let $\varphi:\Omega^*\otimes \mathbb{Q}\longrightarrow R$ be an
arbitrary Hirzebruch genus, i.e. a homomorphism from some bordism
ring to some ring $R$ (usually integers or complex numbers). It
means that $\varphi$ behaves well under disjoint sums and products
of manifolds equipped with some structures (orientation, almost
complex or spin etc.).

Let
${\mathcal{E}}:\Gamma(E_0)\stackrel{D_0}{\longrightarrow}\Gamma(E_1)
\stackrel{D_1}{\longrightarrow} \cdots
\stackrel{D_{q-1}}{\longrightarrow} \Gamma(E_q)$ be an elliptic
complex of differential operators on $X$ whose index $
ind\{\mathcal{E}\}=\sum_{j=0}^q (-1)^j{\rm dim}_{
\mathbb{C}}H^j(\mathcal{E})$ is equal to the genus $\varphi(X)$.
Classical complexes and operators, as the de Rham complex for
oriented manifolds, the Dolbeault complex for complex manifolds
and the Dirac operator for $Spin^c$ manifolds are elliptic. It is
the statement of the Atiyah-Singer Index Theorem \cite{HBJ} that
the indexes of these complexes can be calculated topologically by

\[ind(\mathcal{E})=\Bigl((\sum_{i=0}^q
(-1)^ich(E_i))\prod_{j=1}^n(\frac{x_j}{1-e^{-x_j}}\cdot\frac{1}{1-e^{x_j}})\Bigr)[X],\]
where $x_i$ are the Chern roots of $X$ and $ch(E_i)$ are the Chern
characters of bundles $E_i$. The equivariant genus for manifolds
with finite group action compatible with operators in the given
complex $\mathcal{E}$ on $X$ is defined by \[
\varphi(g,X)=\sum_{j=0}^q(-1)^j{\rm tr}g^*|_{H^j( \mathcal{E})}.\]
Based on the previous formula for the Euler characteristic and
signature we consider the $\varphi-$genus of orbit space as the
averaging sum of equivariant
genera\[\varphi(X/G)=\frac{1}{|G|}\sum_{g\in G}\varphi(g,X).\]
Note the appearance of the cyclic index in calculations of genera
of symmetric powers defined in this way. Namely, \[
\varphi(SP^n(M))=\frac{1}{n!}\sum_{\sigma \in S_n}\varphi
(\sigma,M^n)=\sum_{\alpha_1+2\alpha_2+\cdots
n\alpha_n=n}\frac{\prod_{r=1}^n\varphi(\omega_r,M^r)^{\alpha_r}}{1^{\alpha_1}\cdots
n^{\alpha_n} \alpha_1!\cdots\alpha_n!},\] where the last equality
holds because of the cycle decomposition of permutations. We use
the generating function of cyclic index to obtain the
corresponding generating functions for genera of symmetric powers
\[ \sum_{n=1}^\infty \varphi(SP^n(X))t^n=exp\Bigl(
\sum_{n=1}^\infty \varphi(\omega_r,X^r)\frac{t^r}{r}\Bigr) .\] We
calculate in this way, $\chi_y-$characteristic and elliptic genus
of symmetric powers $SP^n(S_g)$ of complex curves $S_g$.

Recall that $\chi_y-$characteristic is the index of the Dolbeault
complex associated to the complex manifolds $X$ and following
\cite{HBJ} it can be computed by the formula
\[\chi_y(X)=(\prod_{j=1}^n\frac{x_j(1+ye^{-x_j})}{1-e^{-x_j}})[X].\]
For a complex curve $X=S_g$ of genus $g$ it follows from the
equivariant Atiyah-Singer Index Theorem that\[
\chi_y(\omega_r,S_g^r)=x\cdot\Bigl(\frac{1+ye^{-x}}{1-e^{-x}}\cdot
\frac{1+ye^{i\lambda_1-x}}{1-e^{i\lambda_1-x}}\cdots
\frac{1+ye^{i\lambda_{r-1}-x}}{1-e^{i\lambda_{r-1}-x}} \Bigr)[S_g]
.\] We need to find the coefficient of $x$ in the Taylor expansion
of above product, and it turns out to be \[
\chi_y(\omega_r,S_g^r)=(1-g)(1+(-y)^r).\] This gives the following
\cite{VG}
\[ \sum_{n=1}^\infty
\chi_y(SP^n(S_g))t^n=exp\Bigl((1-g)\sum_{r=1}^\infty(1+(-y)^r)\frac{t^r}{r}\Bigr)
=\bigl((1-t)(1+yt)\bigr)^{g-1}.\] In the case of projective line
($g=0$), we know that $SP^n( \mathbb{C}P^1)= \mathbb{C}P^n$, and
our result agrees with the derivative of logarithm for
$\chi_y-$characteristic, which is
\[ g'(t)=\sum_{n=1}^\infty \chi_y(
\mathbb{C}P^n)t^n=\frac{1}{(1-t)(1+yt)}.\] For $y=-1,0,1$ we
obtain the generating functions for the Euler characteristics,
Todd genera and signatures respectively.

The {\em elliptic genus} $Ell(X)$ is introduced by Witten as the
equivariant genus of the natural circle action on the free loop
space $ \mathcal{L}X$ of a closed, oriented manifold $X$. Its
logarithm is given by the elliptic integral \[ g(y)=\int_0^y
(1-2\delta t^2+\varepsilon t^4)^{-1/2}dt.\] The characteristic
power series of the elliptic genus is $Q_{Ell}(x)=\frac{x}{f(x)}$,
where $f(x)$ is the solution of differential equation
$(f')^2=1-2\delta t^2+\varepsilon f^4$. Using its product
expansion we have the formula\[
Ell(X)=\varepsilon^{n/2}\Bigl(\prod_{i=1}^{2n}x_i\frac{1+e^{-x_i}}{1-e^{-x_i}}
\prod_{k=1}^\infty\frac{1+q^ke^{-x_i}}{1-q^ke^{-x_i}}\cdot
\frac{1+q^ke^{x_i}}{1-q^ke^{x_i}}\Bigr)[X].\] Hence, for complex
curve of genus $g$, by the equivariant Atiyah-Singer Index
Theorem, we have
\[Ell(\omega_r,S_g^r)=\varepsilon^{1/4}\Bigl(x\prod_{j=0}^{r-1}
\frac{1+e^{i\lambda_j-x}}{1-e^{i\lambda_j-x}}\prod_{k=1}^\infty
\frac{1+q^ke^{-i\lambda_j+x}}{1-q^ke^{-i\lambda_j+x}}\frac{1+q^ke^{i\lambda_j-x}}
{1-q^ke^{i\lambda_j-x}}\Bigr)[S_g].\] So we are able to determine
the needed coefficient of $x$ in above product as
\[ Ell(\omega_r,S_g^r)= \left\{
\begin{array}{ll}0 & \mbox{ r odd}\\ (2-2g)\varepsilon^{1/4} &
\mbox{ r even} \end{array} \right.\] which gives \[
\sum_{n=1}^\infty
Ell(SP^n(S_g))t^n=\frac{1}{(1-t^2)^{(1-g)\varepsilon^{\frac{1}{4}}}}.
\] Note that in the case $g=0$ it is different from the
logarithm of the elliptic genus.

The above formula was proved in greater  generality in
\cite{Zhou}, \cite{BL}. Here we describe another approach to
orbifold genera, motivated by String theory. There is a definition
of orbifold Euler characteristic for manifolds with group actions
\cite{HH}
\[\chi(X,G)=\sum_{\{g\}}\chi(X^g/C(g)).\]
Generalizing this formula for an arbitrary genus, we define the
corresponding orbifold genera
\[\varphi(X,G)=\sum_{\{g\}}\frac{1}{|C(g)|}\sum_{h\in
C(g)}\varphi(h,X^g).\] In the case of orbifold elliptic genera of
symmetric powers, for a manifold $X$ with elliptic genus
$Ell(X)=\sum_{m,l}c(m,l)y^lq^m$, we have the formula due to
Dijkgraaf, Moore, Verlinde, Verlinde \cite{DiMoVeVe}
\[\sum_{n\geq
0}Ell(X^n,S_n)t^n=\prod_{i=1}^\infty\prod_{l,m}\frac{1}{(1-t^iy^lq^m)^{c(m,l)}}.\]

\end{document}